\documentclass{amsart}

\newcommand\quotient{\backslash}
\newcommand\Aff{\mathord{\operatorname{Aff}(\Lambda\backslash G)}}
\newcommand\Ad{\operatorname{Ad}}
\newcommand{\CU}{C(\mathcal X)}

\def\real{\mathord{\mathbb R}}
\def\rational{\mathord{\mathbb Q}}
\def\integer{\mathord{\mathbb Z}}

\renewcommand{\see}[1]{(see~\ref{#1})}

\theoremstyle{plain}
 \newtheorem{thm}{Theorem}[section]
 \newtheorem{cor}[thm]{Corollary}

 \theoremstyle{definition}
 \newtheorem{defn}{Definition\ignorespaces}
           
 \newtheorem{defnwithref}{Definition}
           
 \newtheorem{ack}{Acknowledgment\ignorespaces}

\newcommand{\thmcopier}[1]{\renewcommand\thethm
 {#1}\addtocounter{thm}{-1}}

\begin{document}

\title[Correction to `Quotients of unipotent
translations']{Correction to \\
 ``Measurable quotients of unipotent translations
 \\ on homogeneous spaces" \\
 $\textnormal{[\frenchspacing\emph{Trans. Amer. Math. Soc. \bf
345} (1994) 577--594]}$}

\author{Dave Witte}

\address{Department of Mathematics, Oklahoma State University, 
 Stillwater, OK 74078}
 \email{dwitte@math.okstate.edu}

\begin{abstract}
 The statements of Main~Theorem~1.1 and Theorem~2.1 of the
author's paper [\emph{Trans.\ Amer.\ Math.\ Soc.}\ {\bf 345}
(1994) 577--594] should assume that $\Gamma $~is discrete and
$G$~is connected. (Cors.~1.3, 5.6, and~5.8 are affected
similarly.) These restrictions can be removed if the conclusions
of the results are weakened to allow for the possible existence of
transitive, proper subgroups of~$G$. In this form, the results
can be extended to the setting where $G$ is a product of real
and $p$-adic Lie groups.
 \end{abstract}


\date{\today}

\thanks{Submitted to the \emph{Transactions of the American
Mathematical Society} in June 1996.}

\maketitle

There are two errors in the proof of Theorem~2.1
of~\cite{Witte-quotients}. To eliminate these mistakes,
 \begin{quote}
 \emph{the statements of Theorem~2.1 and Main~Theorem~1.1 should
assume that\/ $\Gamma$ is discrete and $G$ is connected.}
 \end{quote}
 A few other results are affected: In Cor.~1.3, assume $\Gamma $
is discrete and $U$ is connected. In Cor.~5.6, assume $\Lambda$
is discrete. In Cor.~5.8, assume $\Lambda$ is discrete and $N$~is
connected.

Although it includes the main cases of interest, the restriction
to connected groups and discrete subgroups is not entirely
satisfying. For more general Lie groups, however, there are two
problems in the proof: (1)~it was tacitly assumed that $\Aff$ is
second countable, but this may not be the case if $G/G^\circ$ is
not finitely generated; and (2)~some proper subgroups of~$G$ may
be transitive on $\Gamma \quotient G$, so $X$ may not project onto
all of~$G$ in Step~1 on p.~582. The first mistake can be
eliminated by using the topology of convergence in measure. The
second problem can be eliminated by hypothesis (as above), but it
can be better resolved by weakening the conclusion to account for
the transitive proper subgroups (see Cor.~\ref{unip-quot} below). 

The improved argument that eliminates Problem~(1) applies to any
locally compact group, not just Lie groups \see{general-quot}.
Because M.~Ratner's Classification of Invariant Measures~(1.2) is
now known to be true not just for Lie groups, but also for direct
products of real and $p$-adic Lie groups \cite{Margulis-Tomanov},
\cite{Ratner-S-alg}, this allows us to extend Main Theorem~1.1 to
this more general setting \see{S-unip-quot}.

Before stating our results, let us present a simplified definition
of central double-coset quotients. The original definition
imposed a more complicated restriction on~$K$, because the author
did not realize that noncompact groups are unnecessary (when
$\Gamma \quotient G$ has finite volume).

\begin{defnwithref}[{cf.~\cite[p.~578]{Witte-quotients}}]
 Let $U$ be a subgroup of a locally compact group~$G$, and let
$\Gamma \quotient G$ be a finite-volume homogeneous space of~$G$.
Suppose
 \begin{enumerate}
 \item $\Lambda$ is a closed subgroup of~$G$ containing~$\Gamma
$; and
 \item $K$ is a compact subgroup of $\Aff$ that centralizes~$U$.
 \end{enumerate}
 Then the natural $U$-action on $\Lambda \backslash G/K$ is a
quotient of the $U$-action on $\Gamma \backslash G$. It is called
a \emph{central double-coset quotient}.

We should point out that the topology on $\Aff$ here is the
topology of convergence in measure. The more natural topology of
uniform convergence on compact sets may be used instead if
$\Aff$ is second countable and locally compact in this topology.
Unfortunately, this is not always the case, although it is true
if $G$ is a Lie group whose component group $G/G^\circ$ is
finitely generated.
 \end{defnwithref}

We also present a naive analogue, for non-Lie groups, of the
notion of Lie subgroup of a Lie group.

\begin{defn}
 Let us say that a subgroup~$H$ of a locally compact group~$G$ is
\emph{semi-closed} if, for some~$n$, there is a closed
subgroup~$X$ of $G^n = G \times G \times \cdots \times G$, such
that $H$ is the projection of~$X$ into the first factor of $G^n$.
 \end{defn}

\begin{thm} \label{general-quot}
 Let $U$ be a subgroup of a locally compact, second countable
group~$G$, and let\/ $\Gamma \quotient G$ be a finite-volume
homogeneous space of~$G$, such that every semi-closed subgroup
of~$G$ that contains~$U$, and is transitive on\/ $\Gamma \quotient
G$, is open. Assume that the $U$-action on\/ $\Gamma \quotient G$
is ergodic, and that every ergodic invariant probability measure
for the diagonal $U$-action on\/ $(\Gamma \quotient G) \times
(\Gamma \quotient G)$ is homogeneous for a subgroup of $G \times
G$.  Then each quotient of the $U$-action on\/ $\Gamma \quotient
G$ is isomorphic to a central double-coset quotient of\/ $(\Gamma
\cap H)\quotient H$, for some open subgroup~$H$ of~$G$ that
contains~$U$, and is transitive on\/ $\Gamma \quotient G$.
 \end{thm}

If $G$ is a Lie group, then, by replacing it with a Lie
subgroup of the smallest dimension among those with the property
that they are transitive and contain~$U$, we can arrange for the
theorem's hypothesis on semi-closed subgroups to be fulfilled.
Then, by combining the theorem with Ratner's Theorem~(1.2), we
obtain the following version of Main Theorem~1.1 that does not
assume $\Gamma$ is discrete or that $G$ is connected, but has a
weaker conclusion.

\begin{cor} \label{unip-quot}
 Let $U$ be a nilpotent, unipotent subgroup of a Lie group~$G$,
and let\/ $\Gamma \quotient G$ be a finite-volume homogeneous
space of~$G$. If the $U$-action on\/ $\Gamma \quotient G$ is
ergodic, then each quotient of the $U$-action on\/ $\Gamma
\quotient G$ is isomorphic to a central double-coset quotient
of\/ $(\Gamma \cap H)\quotient H$, for some Lie subgroup~$H$
of~$G$ that contains~$U$, and is transitive on\/ $\Gamma \quotient
G$. \qed
 \end{cor}

If $G$ is a direct product of real Lie groups and $p$-adic Lie
groups (for different~$p$'s), then any closed subgroup of~$G$
contains an open, normal subgroup that is a product of closed
subgroups of the factors of~$G$ \cite[Prop.~1.5,
p.~289]{Ratner-S-alg}. Thus, the dimension of any closed subgroup
(or of any semi-closed subgroup of any semi-closed subgroup
of~\dots of~$G$) can be defined as the sum of the dimensions of
its intersections with the direct factors of~$G$. Thus, by
replacing $G$ with an appropriate subgroup of minimal dimension,
we can again arrange for the theorem's hypothesis on semi-closed
subgroups to be fulfilled. So, because Ratner's Theorem has been
generalized to (many) groups of this type
\cite{Margulis-Tomanov}, \cite{Ratner-S-alg}, we obtain the
following classification of quotients for homogeneous spaces of
these product groups. (See \cite[p.~276]{Ratner-S-alg} for the
definition of \emph{regular}; every algebraic group is regular.)

\begin{cor} \label{S-unip-quot}
 Let $S$ be a finite set of places of\/~$\rational$ and, for each
$p \in S$, let $G_p$ be a regular Lie group over\/~$\rational_p$
\textnormal{(}where $\rational_p =\real$ if $p =
\infty$\textnormal{)}. Let $G$ be a closed subgroup of the direct
product\/ $\times_{p \in S} G_p$.
 Let $U$ be a subgroup of~$G$ that is generated by
one-parameter\/ $\Ad$-unipotent subgroups of\/ $\times_{p \in S}
G_p$, and let\/ $\Gamma$ be a lattice in~$G$. If the $U$-action
on\/ $\Gamma \quotient G$ is ergodic, then each quotient of the
$U$-action on\/ $\Gamma \quotient G$ is isomorphic to a central
double-coset quotient of\/ $(\Gamma \cap H)\quotient H$, for some
semi-closed subgroup~$H$ of~$G$ that contains~$U$, and is
transitive on\/ $\Gamma \quotient G$. \qed
 \end{cor}

\begin{defn}
 Let $\Gamma \quotient G$ be a homogeneous space of a locally
compact group~$G$, and let $A,B \subset G$. We say that a Borel
function $\phi \colon \Gamma \quotient G \to \Gamma \quotient G$
is \emph{affine for~$A$ via~$B$} if, for each $a \in A$, there is
some $b \in B$, such that $\phi (xa) = \phi (x) b$ for all $x \in
\Gamma \quotient G$. Note that if the subgroup~$\Gamma$ contains
no nontrivial, normal subgroup of~$G$, then the element~$b$ is
unique, for a given $a \in A$. Furthermore, if $\phi$ is affine
for~$G$ (via~$G$), then $\phi$ is an affine map.
 \end{defn}

\begin{proof}[\bf Proof of Thm.~\ref{general-quot}.]
 Follow the proof of Thm.~2.1 (with ``locally compact" in the
place of ``Lie") until the definition of~$L$, which is a displayed
equation on p.~582. (There is a typographical error at the start
of line~9 on p.~582: $(S,\nu)$ should be $(T,\nu)$.) 

From the proof of Step~1 (on p.~582), and using the hypothesis
on semi-closed subgroups of~$G$, we see that, for a.e.~$(\Gamma
s,\Gamma t) \in M$, there is a closed subgroup~$X$ of~$G\times G$
such that {\rm (1)}~for almost every $x \in X$, we have $(\Gamma
s, \Gamma t) \cdot x \in M$, {\rm (2)}~the subgroup~$X$ projects
onto an open subgroup of each of the two factors of~$G \times G$,
and {\rm (3)}~the subgroup~$X$ contains~$U$.

Let $\CU$ be the semi-group of measure-preserving Borel maps from
$\Lambda \quotient G$ to $\Lambda \quotient G$ that commute a.e.\
with the action of each $u \in U$, where two such maps are
identified if they agree~a.e.\ \cite[p.~532]{dJ-R}. From the
proofs of Steps~2 and~3, we see that a.e.\ ergodic component of
the projection of~$\rho$ to $\Lambda \quotient G \times \Lambda
\quotient G$ is supported on the graph of some $\phi \in \CU$.
(Furthermore, $\phi$ is affine for some open subgroup of~$G$,
via~$G$.) Then a theorem of Veech \cite[Thm.~1.8.2]{dJ-R} states
that there is a compact subgroup~$K$ of~$\CU$, such that the
original quotient $(T,\nu)$ is isomorphic to $\Lambda \backslash
G/K$ (and, furthermore, we know that a.e.~$\phi \in K$ is affine
for some open subgroup of~$G$, via~$G$). Thus, all that remains
is to show that there is some open subgroup~$H$ of~$G$, such that
$K$ is affine for~$H$ via~$H$, for this implies that $K$ is a
group of affine maps on $(\Gamma \cap H) \quotient H$. 

First, note that it suffices to show that $K$ is affine for some
open subgroup of~$G$. For, if this is the case, we may let 
 $$H = \{ a \in G \mid \textnormal{$K$ is affine for~$a$} \} .$$
 For any $\phi \in K$ and $a \in H$, we know that $\phi$ is
affine for~$a$ via some $b \in G$. For any $\psi \in K$, because
the composition $\psi \phi$ is affine for~$a$, we see that $\psi$
is affine for~$b$. Since $\psi$ is arbitrary, this implies $b \in
H$. So $\phi$ is affine for~$a$ via~$H$, as desired.

Second, note that if two maps are each affine for (perhaps
different) open subgroups of~$G$ (via~$G$), then their composition
is also affine for an open subgroup of~$G$ (via~$G$). Thus,
because $K$ has no proper conull subgroups, we conclude that each
$\phi \in K$ is affine for some open subgroup of~$G$. The problem
is to find a single open subgroup of~$G$ that works uniformly for
all $\phi \in K$.

Let $A_1 \supset A_2 \supset \cdots$ be a chain of compact
neighborhoods of~$e$ in~$G$, whose intersection is~$\{e\}$, and
let $B_1 \subset B_2 \subset \cdots$ be a chain of compact sets
whose interiors exhaust~$G$. For each $N \in \integer^+$, define
 $$K_N = \{ \phi \in K \mid
 \textnormal{$\phi$~and~$\phi^{-1}$ are affine for~$A_N$
via~$B_N$} \} .$$
 We claim that $K_N$ is closed. Suppose $\phi_n \to \phi \in K$.
Given $a \in A_N$, each $\phi_n$ is affine for~$a$ via some $b_n
\in B_N$. Because $B_N$ is compact, we may assume $\{b_n\}$
converges, to some $b \in B$. Then
 $$ \phi(sa) \leftarrow \phi_n(sa) = \phi_n(s) b_n \to \phi(s) b
 \textnormal{\qquad (convergence in measure)} .$$
 So $\phi$ is affine for~$a$ via~$b$. The same argument applies
to~$\phi^{-1}$, so $\phi \in K$, as desired.

Because every open subset of~$G$ contains some~$A_n$, we know
that each $\phi \in K$ is affine for some~$A_n$ (via~$G$). Then
the compactness of~$A_n$ implies that $\phi$ is affine for~$A_n$
via some~$B_m$. Therefore, $\phi \in K_{\max(m,n)}$. Hence, the
Baire Category Theorem implies there is some $N_1$, such that
$K_{N_1}$ has nonempty interior. Let $C = K_{N_1}$; the
compactness of~$K$ (and the fact that $C^{-1} = C$) implies there
is some $r \in \integer^+$, such that $C^r$ is an open subgroup
of~$K$. This implies that $K/C^r$ is finite. Thus, by
enlarging~$N_1$, we may assume that $K_{N_1}$ contains
coset representatives of $K/C^r$. Hence, $C^r = K$.

By modding out the largest normal subgroup of~$G$ that is
contained in~$\Lambda$, let us assume, for simplicity, that
$\Lambda$ contains no nontrivial, normal subgroup of~$G$. Then,
for any $\phi \in C$ and $a \in A_{N_1}$, there is a unique
$b(\phi ,a) \in B_{N_1}$, such that $\phi$ is affine for~$a$
via~$b$. We claim that the map $b \colon C \times A_{N_1} \to
B_{N_1}$ is continuous. Suppose $\phi_n \to\phi$ and $a_n \to a$,
and let $b = b(\phi ,a)$. Assume $b_n = b(\phi_n,a_n)$ converges,
say to $b' \in B_{N_1}$. We have
 $$ \phi(s) b = \phi (sa) \leftarrow
 \phi_n(sa_n) = \phi_n(s) b_n \to \phi(s) b'
 \textnormal{\qquad (convergence in measure)}
 .$$
 So $b'=b$, as desired.

From the conclusion of the preceding paragraph, and the
compactness of~$C$, we see that there is some $N_2$, such that
$C$~is affine for~$A_{N_2}$ via~$A_{N_1}$. Continuing, we
recursively construct a sequence $\{N_i\}$, such that $C$~is
affine for~$A_{N_i}$ via~$A_{N_{i-1}}$. By composition, then $K =
C^r$ is affine for~$A_{N_r}$ (via~$B_{N_1}$).
 \end{proof}

\begin{ack}
 I am pleased to thank Marina Ratner, for
encouraging me to extend my work, and C\'esar Silva, for several
helpful conversations. This research was partially supported by a
grant from the National Science Foundation.
 \end{ack}

\end{document}